\documentclass[11pt]{article}
\usepackage[T1]{fontenc}
\usepackage{lmodern}
\usepackage{amsmath,amssymb,amsthm}
\usepackage[a4paper,hmargin=26mm,top=25mm,bottom=26mm]{geometry}
\usepackage{microtype}
\usepackage[hidelinks]{hyperref}
\hypersetup{pdftitle={Vague rotational symmetry at BCC [100] twist grain boundaries},
  pdfauthor={Author names}}
\newtheorem{theorem}{Theorem}[section]
\newtheorem{lemma}[theorem]{Lemma}
\newtheorem{proposition}[theorem]{Proposition}
\theoremstyle{definition}
\newtheorem{definition}[theorem]{Definition}
\newcommand{\C}{\mathbb C}
\newcommand{\Q}{\mathbb Q}
\newcommand{\R}{\mathbb R}
\newcommand{\Z}{\mathbb Z}
\newcommand{\e}{\varepsilon}
\newcommand{\T}{\mathbb T}
\DeclareMathOperator{\lcm}{lcm}
\DeclareMathOperator{\dist}{dist}
\numberwithin{equation}{section}
\title{Vague rotational symmetry at BCC $[100]$ twist grain boundaries}
\author{Kai Jiang, Suining Xiong\thanks{Email: \texttt{2301110058@pku.edu.cn}}}
\date{}
\begin{document}
\maketitle

\begin{abstract}
In earlier work (\href{https://arxiv.org/abs/2603.21158}{arXiv:2603.21158})
we introduced vague rotational symmetry for proximal coincidence point sets
at grain boundaries. A prescribed rotation must be realizable at every
pair tolerance $\e$ by moving each midpoint of an admitted atomic pair
by at most $\e/2$, half the maximum pair separation. We prove the proposed
classification for the planar model of a body-centred cubic $[100]$ twist
boundary. On $0\leq\theta\leq45^\circ$, the maximal rotational order is
twelve at $30^\circ$, eight at $45^\circ$, and four at every other
misorientation. Vague $n$-fold symmetry forces $\xi_n M=M$ for the Fourier
module $M$, hence $\Z[\xi_n]\subseteq M$ in crystal coordinates.
For every rotation of finite order preserving $M$, we determine the least
uniform displacement for all sufficiently small $\e$, allowing independent
point motions and a free rotation centre. In this regime twelvefold symmetry at $30^\circ$
costs $\e/(2\sqrt3)$, and eightfold symmetry at $45^\circ$ costs zero.
At a twist of approximately $25.53^\circ$, the Fourier module also has
eightfold symmetry, but realizing it in the plane costs $\sqrt2\,\e$.
\end{abstract}

\section{Introduction}

A rotation of an atomic structure is an exact symmetry if it maps the
set of sites onto itself. At a grain boundary, evidence for
such a rotation often begins with a pattern visible in the boundary plane.
For example, Schober and Balluffi resolved orthogonal dislocation grids in
electron micrographs of twist boundaries in gold and measured their
orientations and spacings \cite{SchoberBalluffi}. In simulations of
body-centred cubic (BCC) $[100]$ twist boundaries, twelvefold and eightfold
motifs are visible in the interfacial density at misorientations of
$30^\circ$ and $45^\circ$, respectively \cite[Fig.~9]{Xiong}. To assign
a rotation to the entire interface, one must specify how far sites may
move while still representing the observed pattern.

The distinction already arises in the four stages of bicrystallography
described by Pond and Vlachavas \cite[pp.~100, 113--114]{Pond}.
Superposing the two lattices gives a dichromatic pattern. Incorporating
the lattice complexes and then cutting along the chosen interface plane
gives an ideal bicrystal with a definite symmetry group. The fourth stage
allows the bicrystal to relax. Once individual atoms can move, the ideal
geometry alone does not determine which rotations survive. Even an
arbitrarily small displacement of one site can break an exact rotation
if the other sites in its orbit remain fixed.

The \emph{proximal coincidence point set} (PCPS) of \cite{Xiong} gives
the displacement a length scale supplied by the atoms themselves. A pair
separated by at most $\e$ contributes its midpoint. Each atom of the pair
thus lies within $\e/2$ of that midpoint. Vague symmetry uses this same
radius as the allowed positional uncertainty for each representative
point. For a prescribed rotation we minimize the largest midpoint
displacement over all pointwise adjustments and all choices of rotation
centre. The resulting length measures how far the entire set is from
realizing that rotation. At $30^\circ$, for example, the least displacement
that realizes twelvefold symmetry is exactly $1/\sqrt3$ of the available
radius for all sufficiently small $\e$.

Iwanowski studied uniformly bounded perturbations of finite planar point
sets \cite{Iwanowski}. For a PCPS the uniform bound applies across the
entire plane, and shrinking $\e$ also changes which pairs contribute.
For twist interfaces,
Romeu described an octagonal interface at a $45^\circ$ twist about a
cubic axis equivalent to $[001]$ \cite[Fig.~4]{Romeu}. Ochana and Lifshitz used midpoints of
near coincidences to construct octagonal tilings from two square
lattices at $45^\circ$ and dodecagonal tilings from two triangular
lattices at $30^\circ$ \cite{OchanaLifshitz}. The twelvefold case
considered here uses \emph{two square lattices}. Our earlier work
\cite[Sec.~III.A.3]{Xiong} proposed the $30^\circ/45^\circ$
classification and gave constructions for the exceptional symmetries.
Here we prove the converse. At every other misorientation the quarter-turn
and its powers exhaust the vague rotations, even when a different centre
can be chosen at each tolerance.

\paragraph{The planar model.}
Take the lattice spacing as the unit of length and align the first
crystal with the coordinate axes. Its interfacial square lattice is
$G=\Z[i]\subset\C$ and that of the second crystal is $uG$, where
$u=e^{i\theta}$. For $x,y\in G$ set
\begin{equation}\label{eq:pairs}
 d=x-uy,\qquad p=\frac{x+uy}{2},\qquad
 \Lambda_\e(u)=\left\{\frac{x+uy}{2}: x,y\in G,\ |x-uy|\leq\e\right\}.
\end{equation}
Thus $d$ is the mismatch and $p$ the midpoint. In the three-dimensional
PCPS, this is the reduction when only the central layer of each grain
enters the window. Statements for every $\e$ below concern the planar
family \eqref{eq:pairs}. The square lattice
identifies angles modulo $90^\circ$, and conjugation identifies $\theta$
with $-\theta$. It suffices to take $0\leq\theta\leq45^\circ$.

\begin{definition}\label{def:vague}
Write $\xi_n=e^{2\pi i/n}$. The family $\Lambda_\e(u)$ has \emph{vague
$n$-fold symmetry} if, for every $\e>0$, there is a map
$Y_\e:\Lambda_\e(u)\to\C$ with
\[
 \sup_{p\in\Lambda_\e(u)}|Y_\e(p)-p|\leq\e/2
\]
whose image is invariant under $z\mapsto\xi_n z+b_\e$ for some
$b_\e\in\C$. For $n>1$ the centre is $b_\e/(1-\xi_n)$ and may depend
on $\e$. Order one denotes the identity.
\end{definition}

For a prescribed rotation, define its displacement cost by
\[
 D_n(\theta,\e)=\inf_Y\,\sup_{p\in\Lambda_\e(u)} |Y(p)-p|,
\]
where $Y$ ranges over maps whose images admit the rotation through
$2\pi/n$ about any centre. Definition~\ref{def:vague} requires a map
that realizes the rotation within the allowance $\e/2$ at every tolerance.

\begin{theorem}[Classification]\label{thm:classification}
The vague rotational orders of the planar BCC $[100]$ twist PCPS are
precisely the divisors of
\[
 \begin{cases}
  12,&\theta=30^\circ,\\
  8, &\theta=45^\circ,\\
  4, &\text{all other }\theta\in[0,45^\circ].
 \end{cases}
\]
For each excluded pair $(\theta,n)$ there is an $\e_0>0$ such that no
perturbation in Definition~\ref{def:vague} exists for $0<\e<\e_0$.
\end{theorem}

Section~\ref{sec:fourier} proves the necessary condition
$\xi_n M=M$ for the Fourier module $M=G+uG$.
Section~\ref{sec:angles} uses the displacement bound to determine the
misorientation and proves that bound for arbitrary pointwise
perturbations. Section~\ref{sec:examples} constructs the two exceptional
rotations and gives a twist with an eightfold Fourier module for which
the required displacement exceeds the allowance.

\section{Vague rotations preserve the Fourier module}\label{sec:fourier}

In inner products we identify complex numbers with their real coordinate
columns, so $k^{T}p=\operatorname{Re}(k\bar p)$. Both square lattices are
self-dual for the characters $p\mapsto e^{2\pi i k^{T}p}$. Their reciprocal vectors generate
$M=G+uG=G-uG$. We recover $M$ from the near-coincidence points
themselves. The use of uniformly small phases is related to the
$\varepsilon$-duals of cut-and-project sets studied by Akama and Iizuka
\cite{AkamaIizuka}. In the present model the window shrinks, and the
following direct argument identifies the limiting module.

Write $\|t\|_{\T}=\dist(t,\Z)$ and set
\begin{equation}\label{eq:E}
 \delta_\e(k)=\sup_{p\in\Lambda_\e(u)}\|k^{T}p\|_{\T},
 \qquad E=\{k\in\C:\lim_{\e\downarrow0}\delta_\e(k)=0\}.
\end{equation}

\begin{lemma}[Phase lemma]\label{lem:phase}
For every $u=e^{i\theta}$, $E=M$.
\end{lemma}
\begin{proof}
If $k=g+h$ with $g\in G$ and $h\in uG$, then
$g^{T}x+h^{T}(uy)$ is an integer. Subtracting it from
$k^{T}p$ leaves $(h-g)^{T}d/2$, since $p=x-d/2=uy+d/2$.
Thus $\delta_\e(k)\leq|h-g|\e/2$ and $M\subseteq E$.

Suppose first that $u\notin\Q(i)$. Each $d=x-uy\in M$ determines
$x,y$ uniquely, since any relation $x=uy$ with $y\ne0$ would put
$u$ in $\Q(i)$. Write $p(d)$ for its midpoint. This is an additive,
injective map on $M$. Also $M$ is dense in $\C$. Indeed, the residues
of $ju$ modulo $G$, $j\in\Z$, are distinct in the unit square. Their
differences give arbitrarily short nonzero vectors $v\in M$, and
$vG\subset M$ is a square grid of mesh $|v|$.

For $k\in E$, the homomorphism $d\mapsto k^{T}p(d)\bmod\Z$
from $M$ to $\R/\Z$ is continuous at zero by \eqref{eq:E}. It is
therefore uniformly continuous and extends to the plane. Every
continuous character of $\R^2$ has the form
$d\mapsto h^{T}d\bmod\Z$ for some $h\in\R^2$.
Consequently
\[
 \tfrac12 k^{T}(x+uy)-h^{T}(x-uy)\in\Z
 \qquad(x,y\in G).
\]
Taking $y=0$ gives $k/2-h\in G$ by self-duality, and taking $x=0$
gives $k/2+h\in uG$. Their sum is $k\in M$.

If $u\in\Q(i)$, then $M$ is a lattice and every sufficiently small
window contains only the mismatch $d=0$. For such windows
$\Lambda_\e(u)=G\cap uG$. This intersection is $M^*$, the lattice
dual to $M$, by self-duality of the two square lattices. Equation
\eqref{eq:E} now gives $E=(M^*)^*=M$.
\end{proof}

\begin{theorem}[Spectral necessity]\label{thm:module}
If Definition~\ref{def:vague} holds even along a sequence
$\e_j\downarrow0$, then
\begin{equation}\label{eq:module}
             \xi_n M=M,\qquad \Z[\xi_n]\subseteq M .
\end{equation}
\end{theorem}
\begin{proof}
Put $r=\xi_n$. At a window admitting a symmetric perturbation $Y_\e$,
let $a=\sup_p|Y_\e(p)-p|$. For $p\in\Lambda_\e$ choose standard points
$q_p,q_0\in\Lambda_\e$ whose images are the rotations of
$Y_\e(p),Y_\e(0)$. Then
\[
 |q_p-(rp+b_\e)|\leq 2a,\qquad |q_0-b_\e|\leq 2a,\qquad
 |rp-(q_p-q_0)|\leq 4a .
\]
Taking differences removes the translation even when the centres vary.
Since both $q_p$ and $q_0$ are standard points, the triangle inequality
on $\R/\Z$ and $(r^{-1}k)^{T}p=k^{T}(rp)$ give
\begin{equation}\label{eq:transfer}
        \delta_\e(r^{-1}k)\leq 2\delta_\e(k)+4|k|a .
\end{equation}
For $k\in E$ and $a\leq\e_j/2$ the right-hand side tends to zero along
the sequence. The left-hand side decreases as the window shrinks, so
$r^{-1}k\in E$. Iterating $r^{-1}E\subseteq E$ and using $r^n=1$
gives $rE=E$. Apply Lemma~\ref{lem:phase}. Since $1\in M$, all
integer combinations of powers of $r$ belong to $M$.
\end{proof}

In fact, if $\xi_n M\ne M$ then $\liminf_{\e\downarrow0}D_n(\theta,\e)>0$.
Choose $k\in M$ with $r^{-1}k\notin M$, which exists by finite order.
The left side of \eqref{eq:transfer} has a positive limit, whereas
$\delta_\e(k)\to0$. Taking the infimum over all symmetric perturbations
in \eqref{eq:transfer} gives the assertion.

\paragraph{Possible orders.}
Multiplication by $i$ always preserves $M$. If $rM=M$, the rotations
generated by $i$ and $r$ form the cyclic group of order
$m=\lcm(4,n)$. Thus $\Z[\xi_m]\subseteq M$. By the cyclotomic degree
formula (see also the planar modules in \cite{BaakeGrimm}),
\begin{equation}\label{eq:rank}
          \varphi(m)\leq\operatorname{rank}_{\Z} M\leq4 .
\end{equation}
The multiples of four satisfying \eqref{eq:rank} are $4,8,12$.
Indeed, for $m=2^j l$ with $l$ odd and $j\geq2$ the bound is
$2^{j-1}\varphi(l)\leq4$. In particular orders five and ten would
both require $\varphi(20)=8$. In the commensurate case $M$ has rank
two, so $n$ must divide four. At every misorientation the quarter-turn
is an exact symmetry of the standard set, by $(x,y)\mapsto(ix,iy)$.

\section{Selection of the twist angle}\label{sec:angles}

The displacement allowance selects the angle from each of the two
remaining cyclotomic possibilities. Suppose $n\nmid4$ and
$\xi_n M=M$. By \eqref{eq:rank}, $m=\lcm(4,n)$ is eight or twelve,
and $\Q(i)+u\Q(i)=\Q(\xi_m)$ since the two spaces have dimension four
and the former contains the latter. Write $u'$ for the image of $u$
under the nonidentity automorphism of $\Q(\xi_m)$ fixing $i$. The
displacement test for this twist is
\begin{equation}\label{eq:test}
                      |u+u'|\leq |u-u'| .
\end{equation}
By Proposition~\ref{prop:cost}, \eqref{eq:test} is necessary and
sufficient for realizing the rotation within $\e/2$ at every tolerance.

Since $\xi_m\in M$, write $\xi_m=c+hu$ with $c,h\in G$ and $h\ne0$.
Use a prime for the same automorphism on $\xi_m$. Then
\begin{equation}\label{eq:integer}
 |\xi_m-c|^2=|h|^2\in\Z,\qquad
 \frac{|u+u'|}{|u-u'|}
       =\frac{|\xi_m+\xi_m'-2c|}{|\xi_m-\xi_m'|}.
\end{equation}

\paragraph{The eightfold case.}
Here $\xi_8=(1+i)/\sqrt2$ and $\xi_8'=-\xi_8$. For $c=A+iB$,
\[
             |\xi_8-c|^2=1+A^2+B^2-\sqrt2(A+B).
\]
Integrality forces $B=-A$. The ratio in \eqref{eq:integer} is then
$|c|=\sqrt2|A|$. By \eqref{eq:test} this is at most one, so $c=0$.
It follows that $|h|=1$ and $u=\xi_8/h$. These choices differ only
by a quarter-turn and give the single misorientation $45^\circ$.

\paragraph{The twelvefold case.}
Here $\xi_{12}=(\sqrt3+i)/2$ and $\xi_{12}'=(-\sqrt3+i)/2$. Now
\[
           |\xi_{12}-c|^2=1+A^2+B^2-\sqrt3 A-B .
\]
Integrality forces $A=0$. The ratio in \eqref{eq:integer} becomes
$|1-2B|/\sqrt3$. Inequality \eqref{eq:test} allows just $B=0,1$.
Both give $|h|=1$, and $\xi_{12}-c$ equals either $\xi_{12}$ or its
complex conjugate. This gives the single misorientation $30^\circ$.
The calculation also applies to a proposed threefold or sixfold
rotation. The combination with a quarter-turn in \eqref{eq:rank} was
made on the Fourier module, so no simultaneous symmetry of a chosen
perturbation has been assumed.

\begin{proposition}[Least displacement for a prescribed rotation]\label{prop:cost}
Under the hypotheses of this section, put
\[
                    t=\frac{u+u'}{2(u-u')}.
\]
For every $\e>0$ there is an image invariant under $\xi_n$ with
uniform displacement $|t|\e$. Moreover,
\begin{equation}\label{eq:cost}
       D_n(\theta,\e)=|t|\e
       \qquad \text{if}\qquad
       0<\e<\frac{1}{\,3/2+6|t|\,}.
\end{equation}
In particular the criterion of Definition~\ref{def:vague} holds
exactly when \eqref{eq:test} holds.
\end{proposition}
\begin{proof}
Put $r=\xi_n$ and denote the automorphism above by a prime on every
element of $\Q(\xi_m)$. It has order two. Since $u$ and $r$ lie
outside $\Q(i)$, we have $u'\ne u$ and $r'\ne r$. Also $r'$ is a
primitive $n$th root, hence a power of $r$. Multiplication by $r'$
therefore permutes the admitted labels $H=M\cap\overline B(0,\e)$.

For $d=x-uy$ define
\begin{equation}\label{eq:q}
         q(d)=\frac{u}{u-u'}\,d'=p(d)+td .
\end{equation}
The equality follows by inserting $d'=x-u'y$. The map is injective and
$q(r'd)=r q(d)$, since $(r')'=r$. Its image is an exact representative
for the rotation. Density of $M$ gives $\sup_{d\in H}|d|=\e$, so the
uniform displacement of this construction is $|t|\e$.

It remains to show why independent movements of the points cannot
lower this cost. It suffices to consider a competing map $Y$ of
displacement $a\leq|t|\e$. In this proof write $Y(d)$ for the
perturbed position of $p(d)$. For each $d\in H$, choose a label
$F(d)\in H$ satisfying $Y(F(d))=rY(d)+b$ for the same affine
rotation. We claim that $F(d)=r'd$. The bounds
\[
 |p(F(d))-rp(d)-b|\leq2a,\qquad |p(F(0))-b|\leq2a
\]
follow from the displacement allowance. By \eqref{eq:q},
$|p(r'd)-rp(d)|=|t(r-r')d|\leq2|t|\e$. Consequently the label
$\Delta=F(d)-F(0)-r'd$ satisfies
\[
       |p(\Delta)|\leq 4a+2|t|\e,\qquad |\Delta|\leq3\e .
\]
The two endpoints $p(\Delta)+\Delta/2\in G$ and
$p(\Delta)-\Delta/2\in uG$ now both have length at most
$(6|t|+3/2)\e$. Under \eqref{eq:cost} this is less than the lattice
spacing, so both vanish. Thus $F(d)=F(0)+r'd$. Since $r'H=H$, we
obtain $F(0)+H\subseteq H$. Iterating this inclusion in a bounded
nonempty set gives $F(0)=0$, proving the claim.

Let $c=b/(1-r)$ be the centre chosen by $Y$. Along the orbit of $d$,
for $j=0,\dots,n-1$, we now have
\[
 \begin{aligned}
    Y((r')^j d)&=r^j\bigl(Y(d)-c\bigr)+c,\\
    p((r')^j d)&=r^j q(d)-t(r')^j d .
 \end{aligned}
\]
Subtract these identities, multiply by $(r')^{-j}$, and average over
$j$. The two geometric sums with ratios $r/r'$ and $1/r'$ vanish,
since $r'\ne r,1$. This leaves
\[
     td=\frac1n\sum_{j=0}^{n-1}(r')^{-j}
               \bigl[Y((r')^j d)-p((r')^j d)\bigr].
\]
Every summand has modulus at most $a$. Hence $|td|\leq a$ for
all $d\in H$, giving $|t|\e\leq a$ as required. The admissible
threshold $|t|\leq1/2$ is precisely \eqref{eq:test}.
\end{proof}

\section{Displacement costs for eightfold and twelvefold symmetry}\label{sec:examples}

At $45^\circ$ the rotation $u=\xi_8$ interchanges the two square
lattices, since $u^2=i$. It preserves the distance between endpoints
and acts on midpoints by the same rotation. Thus
\[
                       D_8(45^\circ,\e)=0\qquad(\e>0).
\]
At $30^\circ$ take $u=\xi_{12}$. Here $G+uG=\Z[\xi_{12}]$, since
$u^2=1+iu$ and $u^3=i$. Equation \eqref{eq:q} gives a single
prescription for all near-coincidence pairs,
\begin{equation}\label{eq:explicit}
    (x,uy)\longmapsto \frac{x+uy}{2}+\frac{i}{2\sqrt3}(x-uy).
\end{equation}
Rotate the mismatch vector by a right angle and add
$1/(2\sqrt3)$ of it to the midpoint. The resulting set has exact
twelvefold symmetry for every $\e$. Proposition~\ref{prop:cost}
proves the sharp value
\[
 D_{12}(30^\circ,\e)=\frac{\e}{2\sqrt3}
 \qquad\left(0<\e<\frac{2}{3+2\sqrt3}\right).
\]
This interval is sufficient for optimality. The construction
\eqref{eq:explicit} obeys the allowance $\e/2$ at every window size.
The quarter-turn symmetry and the subgroups of the two rotations just
constructed give all orders listed in Theorem~\ref{thm:classification}.
If an excluded order were attained along windows tending to zero,
Theorem~\ref{thm:module} and the necessity calculation in
Section~\ref{sec:angles} would admit that angle and order. This also
proves the small-window exclusion in the theorem.

An eightfold Fourier module can also occur at a twist where eightfold
symmetry costs more than the allowed displacement. Take
\[
           c=2(1-i),\qquad u=\frac{\xi_8-c}{3}.
\]
Here $|u|=1$ and the reduced misorientation is approximately
$25.5288^\circ$. The identities $\xi_8=c+3u$ and
$\xi_8 u=3i-cu$ show that $\xi_8 M=M$. Nevertheless
\eqref{eq:integer} gives $|t|=|c|/2=\sqrt2$, and therefore
$D_8(\theta,\e)=\sqrt2\,\e$ for all sufficiently small $\e$.
The eightfold rotation of the Fourier module is exact at this twist.
Its least displacement in the plane exceeds the PCPS allowance by a
factor of $2\sqrt2$.


\begin{thebibliography}{10}

\bibitem{AkamaIizuka}
Y. Akama and S. Iizuka,
\href{https://doi.org/10.1080/14786430701373698}{Cut-and-project sets and
their $\epsilon$-duals}, \emph{Philos. Mag.} \textbf{87} (2007), 2847--2854.

\bibitem{BaakeGrimm}
M. Baake and U. Grimm,
\href{https://doi.org/10.1524/zkri.219.2.72.26322}{Bravais colourings of
planar modules with $N$-fold symmetry}, \emph{Z. Kristallogr.}
\textbf{219} (2004), 72--80.

\bibitem{Iwanowski}
S. Iwanowski,
\href{https://doi.org/10.1016/0304-3975(91)90389-J}{Testing approximate
symmetry in the plane is NP-hard}, \emph{Theoret. Comput. Sci.}
\textbf{80} (1991), 227--262.

\bibitem{OchanaLifshitz}
M. Ochana and R. Lifshitz,
\href{https://doi.org/10.1103/xjzf-xxns}{Tiling by near coincidence},
\emph{Phys. Rev. Research} \textbf{8} (2026), 023330.

\bibitem{Pond}
R. C. Pond and D. S. Vlachavas,
\href{https://doi.org/10.1098/rspa.1983.0028}{Bicrystallography},
\emph{Proc. R. Soc. Lond. A} \textbf{386} (1983), 95--143.

\bibitem{Romeu}
D. Romeu,
\href{https://doi.org/10.1103/PhysRevB.67.024202}{Interfaces and quasicrystals
as competing crystal lattices: Towards a crystallographic theory of interfaces},
\emph{Phys. Rev. B} \textbf{67} (2003), 024202.

\bibitem{SchoberBalluffi}
T. Schober and R. W. Balluffi,
\href{https://doi.org/10.1080/14786437008238400}{Quantitative observation of
misfit dislocation arrays in low and high angle twist grain boundaries},
\emph{Philos. Mag.} \textbf{21} (1970), 109--123.

\bibitem{Xiong}
S. Xiong, W. Zou, P. Zhang, and K. Jiang,
Framework for quasiperiodic interfaces: proximal coincidence point set and
computation, \href{https://arxiv.org/abs/2603.21158}{arXiv:2603.21158} (2026).

\end{thebibliography}
\end{document}